\documentclass[12pt]{amsart}
\usepackage[margin=1in]{geometry}
\usepackage{mathptmx}
\usepackage{mathtools}
\usepackage{amsthm}
\usepackage{amssymb}
\usepackage[alphabetic]{amsrefs}
\usepackage{graphicx}
\usepackage{tikz}
	\usetikzlibrary{decorations.pathreplacing}
	\usetikzlibrary{patterns}
\usepackage{caption}

\newcommand{\ra}{\rightarrow}
\newcommand{\pr}{\prime}
\newcommand{\de}{\partial}

\newcommand{\C}{\mathbb{C}}

\newcommand{\R}{\mathbb{R}}
\newcommand{\Z}{\mathbb{Z}}

\newcommand{\lbar}[1]{\overline{#1}}

\newtheorem{thm}{Theorem}

\theoremstyle{definition}

\newtheorem*{defin*}{Definition}
\theoremstyle{remark}

\begin{document}

\title{PC4 at Age 40}

\author{Michael Freedman}
\address{\hskip-\parindent
	Michael Freedman \\
	Microsoft Research, Station Q, and Department of Mathematics \\
	University of California, Santa Barbara \\
	Santa Barbara, CA 93106 \\
}

\maketitle

This article is not a proof of the Poincar\'{e} conjecture but a discussion of the proof, its context, and some of the people who played a prominent role. It is a personal, anecdotal account. There may be omission or transpositions as these recollections are 40 years old and not supported by contemporaneous notes, but memories feel surprisingly fresh. I have not looked up old papers to check details of statements; this article is merely a download from my current mental state.\footnote{This is a written version for a CMSA lectured delivered at Harvard on September 28, 2020.}
\bigskip

Poincar\'{e} liked to argue in many ways: analytically, combinatorially, and topologically. He seemed averse to even fixing a definition for the term ``manifold,'' so one should not impose modern notions like TOP, PL, and DIFF on his famous conjecture. He himself seems to have thought little about his own conjecture, since it asked if simply connected manifolds were spheres; clearly he meant to also specify that the homology should also vanish (except in the highest and lowest dimensions). The modern statement, now known to be true in every dimension is: a closed topological manifold $M$ that is homotopy equivalent to a sphere is homeomorphic to that sphere. One cannot replace homeomorphism with diffeomorphism throughout the statement because of examples like Milnor's exotic 7-spheres. In all dimensions you still experience most of the excitement if you presume from the start that the manifold $M^n$ has a smooth structure and the goal is to prove it homeomorphic to $S^n$. We will take this perspective.\footnote{In dimension $>3$ work of Lashoff \cite{lashof70} smooths a connected, contractible, topological $M$ and this can be used to bridge to the strongest statement, but the cost is one must use a technically more difficult (proper) version of the work of Smale discussed below. In dimensions $<4$ there are canonical smooth structures \cites{moise52a,moise52b}.}

When the dimension $n$ is 0 or 1 there is not much to prove. The $n=2$ result is a special case of uniformization. The first modern work is Smale's \cite{smale56} in which he proved PC$n$, $n \geq 5$. I will need to say something about his proof since the case of PC4 starts the same way but encounters a special problem when $n=4$. Solving this problem brings in the topological category in a way not present in Smale's work. The reader might think Smale must also have labored in the topological category since the conclusion is only homeomorphism not diffeomorphism. But 99\% of Smale's work is in a smooth setting; his great achievement, the $h$-cobordism, theorem is smooth category. The step from the $h$-cobordism theorem to PC$n$ is small and in some dimensions involves a gluing along an $S^{n-1}$ which might not be extended to a diffeomorphism over $D^n$. So PC4 amounts to Smale's outline with a topological twist. But here the tail wags the dog. When you delve into this detail the twist expands to fill your entire field of view. The final case (historically), dimension 3, was proved by Perelman \cites{perelman03a,perelman03b} using Hamilton's theory of Ricci flow. It is entirely different in outline, more like Beethoven's 9th than a conventional proof, and still stands as the greatest accomplishment of 21st century mathematics.

I worked towards PC4 from 1974 to 1981, roughly half time. I had many mental pictures but no notes when it came to writing the proof. I realized I had no letters (in my mind) for any of the spaces, maps, or relations. Apparently as a youngster I had not learned yet to think in or with symbols. When written, the proof was hard for others to understand; when I tried to explain it, I did not know where to begin. To this day I regret that in my lectures at Berkeley with Smale in the audience, I never even mentioned that it was his proof, with a twist, that I was presenting. From the perspective of a youth in 1981, 1959 might as well have been 1859, I did not feel the historical connections. Now they seem obvious and I will try to capture them here.

Three streams of thought enter the proof and I will personify them with the names of three mathematicians: Steve Smale, Andrew Casson, and Bob Edwards. Of course, each represents a field and there are other names as well, but in all three cases the power of their individual ideas is so strong and so determinative of the form of the final proof that little is lost by identifying them with their fields of work. To put a name to their work (as relevant here):

\begin{itemize}
    \item Smale: Classical smooth topology and dynamics.
    \item Casson: How to get started in 4D: ``finger moves'' simplify $\pi_1$, leading eventually to ``Casson handles.''
    \item Edwards: Manifold factors, decomposition spaces, shrinking.
\end{itemize}

I will explain, in order, the input from each of these three streams and recall the occasional anecdote.

By the time (1959) of Smale’s paper Morse theory had already moved on to infinite dimensions. People knew the Morse inequalities and that seemed to be regarded as a satisfactory link between critical points and homology of manifolds. Smale took a much closer look. Let us recall his set up. Henceforth, all manifolds are assumed to be compact and smooth. Let $W^{n+1}$ be an $(n+1)$-dimensional manifold with two boundary components, $M_0$ and $M_1$. Let us assume all three fundamental groups are trivial, and the relative homology vanishes: $H_\ast(W,M_0;\Z) = 0$. (Of course by duality, the relative homology groups also vanish if we replace $M_0$ with $M_1$.) With these hypotheses $W$ is called a (simply connected) $h$-cobordism.

\begin{thm}[Smale's $h$-cobordism theorem (hCT)]
    A simply connected $h$-cobordism $W^{n+1}$ is diffeomorphic to a product $M_0 \times I$, provided $n > 4$.
\end{thm}

The idea of the proof is to relentlessly match the algebra of the Morse complex of a smooth Morse function, $f:W \ra [0,1]$ to the geometry, or as Smale would say, the dynamics of the gradient ($\nabla f$) flow. The Morse complex has critical points of index $k$ as its chain group generators in dimension $k$, and the differential of the complex records in matrix form the algebraic number $\#_{i,j}$ of times the boundary of the $i$th $k$-cell wraps over the $j$th $(k-1)$-cell. The Morse cancelation lemma states that if the geometric count $\#_{i,j} = \pm 1$, the two critical points can be cancelled, and the Morse functions simplified. Smale realized that the in the simply connected case, the vanishing of the relative homology groups was the only obstruction to performing a series of deformations to $f$ (called ``handle slides,'' or in the algebra, ``row operations'') that lead to cancelation of \emph{all} critical points. Once there are no critical points left, the gradient flow lines define the desired product structure.

Matching the geometry to the algebra consists of turning the algebraic information that the signed sum of intersection points $\#_{i,j} = 1 \in \Z$ into the stronger geometric information that the $i$th descending sphere truly meets the $j$th ascending sphere in one transverse point, not for example three points in a $+, -, +$ pattern. To accomplish this he employed Whitney’s trick \cite{whitney44} for removing pairs of oppositely signed double points. These double points live in an $n$-dimensional cross-section $L^n$ of $W$, often called the middle level, and Whitney’s method requires $n>4$. The dimension restriction arises because a 2 dimensional disk, the ``Whitney disk,'' guides the cancellation of the extra $(+,-)$ pair and we need that disk to be embedded; an easy thing when $\dim(L) = n > 4$, and a difficult thing otherwise. The proof of PC4 follows Smale up to this point, from there forward it is all about how to locate the Whitney disks necessary to make $\{A\}$ and $\{D\}$, the ascending and descending 2-spheres in the middle level $L$, intersect each other \emph{geometrically} in single transverse points, as the algebra would indicate. So Smale has gotten us started but we are left with the problem of simplifying, say, a canonical picture in the middle level $L^4$ of two embedded 2-spheres $A$ and $D$ which meet in three points $+,-,+$.

It is here where Casson’s work begins. Before turning to Casson, this is a good place to explain how the hCT implies various PCn. There are actually two overlapping techniques (both due to Smale). Suppose $\Sigma^n$ is a homotopy $n$ sphere. We can cut out two disjoint balls from $\Sigma^n$ to obtain an $n$ dimensional $h$-cobordism $W^n$. If $n$ is at least 6 the hCT recognizes $W$ as diffeomorphic to $S^{n-1} \times I$. $\Sigma^n$ is now obtained by gluing back the two balls, but the gluing diffeomorphisms reflect the complexity of the product structure output of the hCT and we don’t know that the result is diffeomorphic to $S^n$. We see from this picture that $\Sigma$ is two balls glued together by a homeomorphism of their boundaries. It is elementary (via Alexander's coning trick) that the union is homeomorphic to $S^n$. To handle the lowest dimensional case treated by Smale, $\Sigma^5$, a different strategy is required. Some algebraic topology and a bit of surgery theory is used to construct a 6D $h$-cobordism $W^6$ between $\Sigma^5$ and $S^5$ and the hCT is applied to $W^6$. Interestingly, this proves more, that $\Sigma^5$ is diffeomorphic, not just homeomorphic, to $S^5$.  This approach using surgery only succeeds in finding $W$ in some dimensions but when it works the conclusion is stronger. In our case $n=4$, the surgery approach actually does work and we find a smooth $h$-cobordism $W^5$ from $\Sigma^4$ to $S^4$. The reason that our proof concludes only a homeomorphism at the end of the day is that there is a terrific struggle involving lots of limits and point set topology to find the Whitney disks and they end up being merely topological category objects. The Whitney disks we find spoil the smooth category context of the rest of the proof. It is an astonishing historical coincidence that within 6 months of PC4, Simon Donaldson had figured out enough about which four dimensional homotopy types do not contain smooth manifolds to know for sure that the Whitney disks I built cannot generally be smoothed. These Whitney disks are all that stand in the way of realizing the quadratic form $E_8 \oplus E_8$ as the intersection form of a close smooth 4-manifold, which is excluded by Donaldson’s ``Theorem C'' \cite{donaldson86}. So the excursion into the topological category is necessary if one speaks generally in terms of 4-manifold constructions. Of course there may be some as yet unimagined workaround---or entirely new method---satisfactory to the study of a smooth homotopy sphere $\Sigma^4$. I think it safe to say that the greatest open problem in topology (``The last man standing'') is the smooth category PC4: Is every homotopy 4-sphere diffeomorphic to $S^4$?

Now to Casson. Casson’s first great insight was that there was some hope for Whitney disks in 4D. He realized that the problem of finding an embedded Whitney disk is qualitatively different from finding a slice disk for a classical knot in $S^3$. That prototypical 2-disk embedding problem had been studied starting in the 50s by Fox and Milnor \cite{fox66}: Given a knot $K$ in $S^3$ does it bound a smoothly embedded disk in the 4-ball $B^4$? This problem is infinitely rich, and generally has a negative answer. But Casson realized that Whitney disks in all applications would have homological duals, so by standing farther back, and taking more of the manifold topology into account, the embedding problem might turn out to be less obstructed than in knot theory. This hope bore fruit.  It is easy to show that these three spaces are simply connected: The middle level $L$,  and the complements $L \setminus A$, and $L \setminus D$.  To look for a Whitney disk $D$ in $L$ meeting $A \cup D$ only where it is supposed to, long the boundary of $D$, we would also like $L \setminus (A \cup D)$ to be simply connected. The duality Casson observed implies that $\pi_1(L \setminus (A \cup D))$ is perfect. He realized that it might help to make things worse before trying to make them better. If one takes ones finger and pushes a generic patch of $A$ along an arc and through a generic point on $D$ this ``finger move'' will create two new $+,-$ points of intersection. Each such point will microscopically look like the intersection of two complex lines in $\C^2$ and linking these complex lines is a ``Clifford torus.'' The algebraic point is that the top cell of these Clifford tori are relations saying the two loops linking the two sheets (complex lines in the model) commute in $\pi_1$ of their joint complement. If you take a perfect group and start forcing pairs of elements to commute pretty soon you have a trivial group. In this way Casson got off the ground, by doing enough finger moves between $A$ and $D$ to kill $\pi_1(L \setminus (A \cup D))$. This enabled him to at least immerse the Whitney disk $W$ (with the required normal framing data) he was looking for \cite{casson86}.

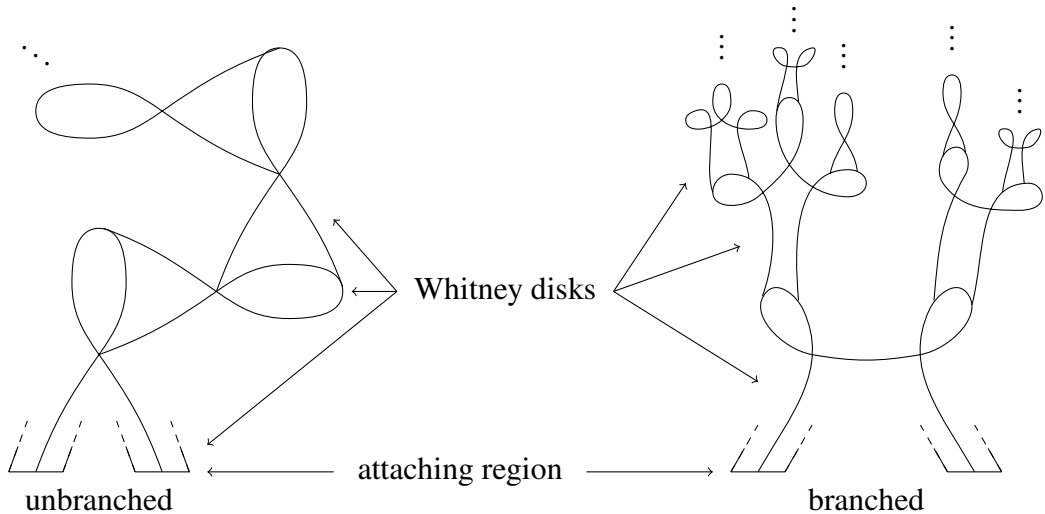
\begin{figure}[ht]
    \centering
    \begin{tikzpicture}[scale=1.2]
        \draw (-0.7,0) to [out=70,in=-125] (0,1.3) to[out=55,in=-90] (0.3,2.1) to[out=90,in=0] (0,2.7) to[out=180,in=90] (-0.3,2.1) to[out=-90,in=125] (0,1.3) to[out=-55,in=110] (0.7,0);
        \draw[yshift=2cm,rotate=-90] (-0.7,0.05) to [out=70,in=-125] (0,1.3) to[out=55,in=-90] (0.3,2.1) to[out=90,in=0] (0,2.7) to[out=180,in=90] (-0.3,2.1) to[out=-90,in=125] (0,1.3) to[out=-55,in=110] (0.7,0);
        \draw[xshift=2cm,yshift=2cm] (-0.7,0) to [out=70,in=-125] (0,1.3) to[out=55,in=-90] (0.3,2.1) to[out=90,in=0] (0,2.7) to[out=180,in=90] (-0.3,2.1) to[out=-90,in=125] (0,1.3) to[out=-55,in=110] (0.7,0.05);
        \draw[yshift=4cm,xshift=2cm,rotate=90] (-0.7,0) to [out=70,in=-125] (0,1.3) to[out=55,in=-90] (0.3,2.1) to[out=90,in=0] (0,2.7) to[out=180,in=90] (-0.3,2.1) to[out=-90,in=125] (0,1.3) to[out=-55,in=110] (0.7,0);
        \node at (-0.7,4.7) {$\ddots$};
        \draw (-1,0) -- (-0.4,0) -- ++(70:0.3);
        \draw (-1,0) -- ++(70:0.3);
        \draw[dashed] (-1,0) -- ++(70:0.6);
        \draw[dashed] (-0.4,0) -- ++(70:0.6);
        \draw (0.4,0) -- (1,0) -- ++(110:0.3);
        \draw (0.4,0) -- ++(110:0.3);
        \draw[dashed] (0.4,0) -- ++(110:0.6);
        \draw[dashed] (1,0) -- ++(110:0.6);
        
        \node at (4,0) {attaching region};
        \draw[<-] (1.2,0) -- (2.6,0);
        \draw[->] (5.4,0) -- (6.8,0);
        
        \draw (7,0) -- (7.6,0) -- ++(60:0.3);
        \draw (7,0) -- ++(60:0.3);
        \draw[dashed] (7,0) -- ++(60:0.6);
        \draw[dashed] (7.6,0) -- ++(60:0.6);
        \draw (7.3,0) to[out=60,in=-65] (7.8,1.8) to[out=115,in=45] (7.4,2) to[out=225,in=170] (7.9,1.3) to[out=-10,in=190] (9.1,1.3) to[out=10,in=-45] (9.6,2) to[out=135,in=65] (9.2,1.8) to[out=-115,in=120] (9.7,0);
        \draw (10,0) -- (9.4,0) -- ++(120:0.3);
        \draw (10,0) -- ++(120:0.3);
        \draw[dashed] (10,0) -- ++(120:0.6);
        \draw[dashed] (9.4,0) -- ++(120:0.6);
        
        \draw (7.35,1.95) to[out=70,in=0] (7,3.3) to[out=180,in=90] (6.8,3.1) to[out=-90,in=220] (7.3,3.05) to[out=40,in=-90] (7.8,3.9) to[out=90,in=0] (7.65,4.15) to[out=180,in=90] (7.5,3.9) to[out=-90,in=180] (8.3,3.05) to[out=0,in=-90] (8.5,3.2) to[out=90,in=0] (8.3,3.35) to[out=180,in=95] (7.75,1.9);
        \draw (9.25,1.9) to[out=85,in=240] (9.6,3.3) to[out=60,in=0] (9.5,3.6) to[out=180,in=125] (9.35,3.2) to[out=-55,in=180] (10.2,2.9) to[out=0,in=-90] (10.45,3.05) to[out=90,in=0] (10.3,3.2) to[out=180,in=75] (9.67,1.85);
        \draw (6.8,3.1) to[out=110,in=20] (6.6,4) to[out=200,in=90] (6.5,3.9) to[out=-90,in=225] (6.9,3.9) to[out=45,in=-90] (7,4.1) to[out=90,in=0] (6.9,4.3) to[out=180,in=90] (6.8,4.1) to[out=-90,in=180] (7.2,3.8) to[out=0,in=-90] (7.4,3.9) to[out=90,in=10] (7.2,4) to[out=200,in=95] (7.2,3.25);
        \node at (6.9,4.8) {$\vdots$};
        \draw (7.51,4) to[out=85,in=10] (7.5,4.7) to[out=200,in=170] (7.6,4.5) to[out=-10,in=190] (7.8,4.5) to[out=10,in=10] (7.9,4.7) to[out=190,in=95] (7.75,4.1);
        \node at (7.7,5.1) {$\vdots$};
        \draw (8.1,3.3) to[out=80,in=-90] (8.35,4) to[out=90,in=0] (8.25,4.2) to[out=180,in=90] (8.15,4) to[out=-90,in=95] (8.4,3.33);
        \node at (8.25,4.7) {$\vdots$};
        
        \draw[xshift=1.2cm,yshift=0.2cm] (8.15,3.3) to[out=80,in=-90] (8.35,4) to[out=90,in=0] (8.25,4.2) to[out=180,in=90] (8.15,4) to[out=-90,in=95] (8.4,3.33);
        \node at (9.45,4.9) {$\vdots$};
        \draw[xshift=2.5cm,yshift=-0.9cm] (7.51,4) to[out=85,in=10] (7.5,4.7) to[out=200,in=170] (7.6,4.5) to[out=-10,in=190] (7.8,4.5) to[out=10,in=10] (7.9,4.7) to[out=190,in=95] (7.75,4.1);
        \node at (10.2,4.2) {$\vdots$};
        
        \node at (0,-0.3) {unbranched};
        \node at (8.5,-0.3) {branched};
        \node at (4.5,2) {Whitney disks};
        \draw[->] (3.3,2) -- (1.2,0.3);
        \draw[->] (3.3,2) -- (2.8,2);
        \draw[->] (3.3,2) -- (2.6,2.8);
        \draw[->] (5.7,2) -- (7.3,1);
        \draw[->] (5.7,2) -- (7.1,2.5);
        \draw[->] (5.7,2) -- (6.5,3.2);
    \end{tikzpicture}
    \caption{Casson handles.}
    \label{fig:cassonhandle}
\end{figure}

From here the finger moves simply explode in number. Fearlessly, Casson decides to cure the problem of the Whitney disk double points by capping them off with subsidiary Whitney disks in a hierarchy later called a Casson tower. When taken to its infinite limit, the open, tapered, regular neighborhood is called a ``Casson handle.'' I draw a dimensionally reduced schematic of a Casson handle in Figure \ref{fig:cassonhandle}.

Casson (1974) showed that his Casson handles $C$ had the correct proper homotopy type to serve as neighborhoods of embedded Whitney disks, $C \overset{P}{\simeq} (D^2 \times \R^2, \de D^2 \times \R^2)$. My contribution seven years later was to replace proper homotopy equivalence---which was not strong enough for most geometric applications---with homeomorphism. As I mentioned, I spent those seven years staring at his handles about half time. I think this was the correct fraction. I was young and needed to do some projects likely to work in order to get tenure, raise a family, and preserve sanity. Those years I was either climbing rocks, writing pleasant papers, or struggling with Casson handles. Not much else. My mother had taught me that everything must be fun, ``fun first, pleasure second'' as she put it. My father who was always at work on something (and always enjoyed it) told me that there is ``your work and, then, the world’s work.'' The idea of course was to do both. My papers were the ``world’s work.'' I enjoyed them greatly and would surely have lost confidence in myself if I was not publishing regularly. Casson handles were ``my work.'' The climbing sheer pleasure. I had a reputation among my friends for ``not thinking on lead,'' I just would tear ahead and see what happened. That is the way I liked it: When I did math I tried to think (a bit); when I climbed I just let it rip.  Many of my friends asked me variations of, ``Since Casson is obviously so much smarter than you, why was it that he didn’t analyze his own handles.'' The answer, I think, is ``opportunity costs''; he was so full of brilliant ideas: secondary obstructions to knot cobordism (with Gordon), the Casson invariant, what became ``weak Heegaard reduction'' (also with Gordon), that he was simply too busy.

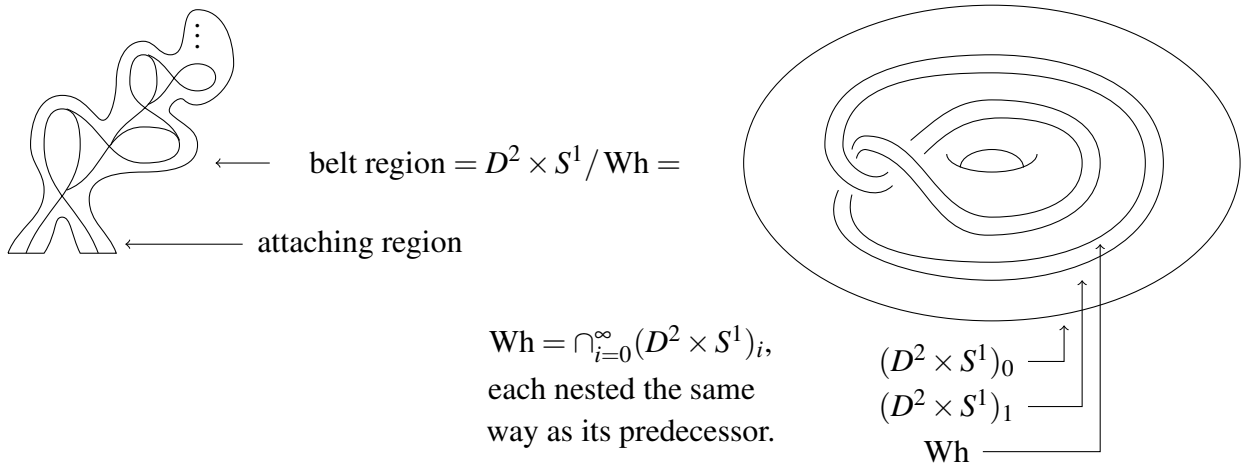
\begin{figure}[ht]
    \centering
    \begin{tikzpicture}[scale=1.2]
        \draw (-3.7,0) to[out=70,in=-90] (-3.1,1.2) to[out=90,in=0] (-3.3,1.6) to[out=180,in=90] (-3.5,1.2) to[out=-90,in=110] (-2.9,0);
        \draw (-3.25,1.6) to[out=-20,in=180] (-2.4,1) to[out=0,in=-90] (-2,1.2) to[out=90,in=0] (-2.4,1.4) to[out=180,in=20] (-3.25,0.7);
        \draw (-2.76,1.2) to[out=70,in=-90] (-2.25,1.9) to[out=90,in=0] (-2.4,2.2) to[out=180,in=90] (-2.55,1.9) to[out=-90,in=110] (-2,1.25);
        \draw (-2.35,2.2) to[out=-20,in=180] (-1.8,1.8) to[out=0,in=-90] (-1.6,1.95) to[out=90,in=0] (-1.8,2.1) to[out=180,in=20] (-2.35,1.65);
        \node at (-1.8,2.5) {$\vdots$};
        \draw (-3.9,0) -- (-3.5,0) .. controls (-3.35,0.23) and (-3.4,0.4) .. (-3.3,0.4) .. controls(-3.2,0.4) and (-3.2,0.2) .. (-3.1,0) -- (-2.7,0) to[out=110,in=-90] (-3,0.6) to[out=90,in=180] (-2.4,0.9) to[out=0,in=-90] (-1.8,1.2) to[out=90,in=-45] (-2.1,1.55) to[out=135,in=135] (-1.95,1.7) to[out=-45,in=-90] (-1.4,2.1) to[out=90,in=0] (-1.8,2.7) to[out=180,in=0] (-2.1,2.2) to[out=180,in=0] (-2.4,2.35) to[out=180,in=90] (-2.7,1.9) to[out=-90,in=-45] (-2.9,1.5) to[out=135,in=0] (-3.3,1.75) to[out=180,in=90] (-3.65,1.2) to[out=-90,in=90] (-3.5,0.7) to[out=-90,in=70] (-3.9,0);
        
        \node at (1.5,1) {belt region $= D^2 \times S^1 \slash \operatorname{Wh} =$};
        \draw[->] (-1,1) -- (-1.6,1);
        
        \draw (7,1) ellipse (2.75 and 1.75);
        \draw (7.5,1.1) arc (0:-180:0.5 and 0.2);
        \draw (7.35,0.95) arc (0:180:0.35 and 0.2);
        \draw (5.45,1.15) to[out=80,in=180] (5.6,1.3) to[out=0,in=180] (7,0.4) to[out=0,in=-90] (8,1) to[out=90,in=0] (7.1,1.5) to[out=180,in=45] (6.25,1.1);
        \draw (5.5,1.05) to[out=80,in=180] (5.6,1.15) to[out=-10,in=180] (6.95,0.2) to[out=0,in=-90] (8.2,1) to[out=90,in=0] (7,1.7) to[out=180,in=40] (6.1,1.25);
        \draw (5.85,0.9) to[out=225,in=-90] (5.35,1.2) to[out=90,in=180] (7,2) to[out=0,in=90] (8.7,1) to[out=-90,in=0] (7,-0.1) to[out=180,in=255] (5.45,0.65);
        \draw (5.9,0.75) to[out=225,in=-90] (5.15,1.2) to[out=90,in=180] (7,2.2) to[out=0,in=90] (8.9,1) to[out=-90,in=0] (7,-0.3) to[out=180,in=245] (5.3,0.7);
        
        \node at (0,0.1) {attaching region};
        \draw[->] (-1.3,0.1) -- (-2.6,0.1);
        
        \node at (3,-1) {$\operatorname{Wh} = \cap_{i=0}^\infty (D^2 \times S^1)_i$,};
        \node at (2.9,-1.5) {each nested the same};
        \node at (3,-2) {way as its predecessor.};
        
        \node at (6.5,-1.2) {$(D^2 \times S^1)_0$};
        \draw[->] (7.4,-1.2) -- (7.8,-1.2) -- (7.8,-0.8);
        \node at (6.5,-1.7) {$(D^2 \times S^1)_1$};
        \draw[->] (7.4,-1.7) -- (8,-1.7) -- (8,-0.3);
        \node at (6.5,-2.2) {$\operatorname{Wh}$};
        \draw[->] (6.9,-2.2) -- (8.2,-2.2) -- (8.2,0.1);
    \end{tikzpicture}
    \caption{$(D^2 \times S^1 \slash \operatorname{Wh}) \times \R$ is homeomorphic to $D^2 \times S^1 \times \R$.}
    \label{fig:beltregion}
\end{figure}

Let us turn to the third stream, Bob Edwards and the world of decomposition space theory, the Bing shrinking criterion (BSC), and manifold factors. In 1980, Bob was the undisputed heavyweight champ in what might variously be called Texas-topology, Bing-topology, the R.L.Moore-R.H.Bing school.  It is an area at the wild end of manifold topology, with roots and emphasis in low dimensions. The Alexander horned sphere is what you should picture. It drew on point set topology but was a very pictorial, concrete, and vivid undertaking. It also had an enjoyably sportive feel that allows me to call Bob ``Champ'' whereas it would never occur to me to call Hiranaka or Deligne the ``Champ'' of algebraic geometry. During those seven years, I had discovered something, now called re-imbedding. If you have a Casson tower of some critical minimum height, then 7, later reduced to 5 (see \cites{cha16,gs84}), you can leave the lowest stage alone and build a new extension inside a neighborhood of the original stage construction which has any, arbitrarily large, number of stages. This quickly allows geometric control to be added to Casson handles. The original construction just kind of flops around, the higher stages are in no sense smaller than the lower ones. But with re-imbedding it is only a small step to build Casson towers in which the higher stages converge beautifully to a Cantor set. Now if one takes a tapered regular neighborhood and completes it with this limiting Cantor set one obtains a compactified $\lbar{C}$ in any application where formerly a Casson handle $C$ was present. A handle has its boundary divided into  the attaching region and the co-attaching or ``belt'' region. See Figure \ref{fig:beltregion}. In early fall of 1980 I used the Kirby calculus to draw an exact picture of the belt region for the standard unbranching $\lbar{C}$. For the everyday closed 2-handle $(D^2 \times D^2, S^1 \times D^2)$, the belt region is $D^2 \times S^1$, a solid torus. What I drew belt ($\lbar{C}$) was what I later learned was a decomposition space $D^2 \times S^1 \slash \operatorname{Wh}$, the solid torus with a compactum called the Whitehead continuum crushed to a point (and endowed with the quotient topology). A few weeks later I was sitting at a pizza restaurant after a meeting of the recurring Southern California Topology meeting with Bob Edwards. Apropos of nothing in particular, he was putting pen to napkin to explain the ``Andrews-Rubin'' \cite{andrews65} shrinking argument which proves that the space $D^2 \times S^1 \slash \operatorname{Wh}$ is a manifold factor. Even the concept to me was astonishing: $D^2 \times S^1 \slash \operatorname{Wh}$ is most certainly not a manifold but it turns out that its product with the real line $(D^2 \times S^1 \slash \operatorname{Wh}) \times \R$ is homeomorphic to $D^2 \times S^1 \times \R$.  Edwards understood this phenomenon well as he had just completed his proof that the double suspension of the homology sphere bounding the Mazur manifold was $S^5$. The single suspension is of course not a manifold. My excitement was overwhelming. The frontier of my controlled Casson handle $\lbar{C}$, while not a manifold, was close---if it turned out to have a product collar neighborhood that would, by Andrews-Rubin, be a manifold and would harbor the Whitney disk I had been seeking for six years. In the end, the proof proceeded somewhat differently but the fact that the belt region of ``frontier'' of the controlled Casson handle was a manifold factor convinced me that PC4 was likely true.

To appreciate why $D^2 \times S^1 \slash \operatorname{Wh}$ is a manifold factor we introduce the Bing Shrinking Criterion (BSC), as exposited by Edwards \cite{edwards80}.

Let $X$ and $Y$ be compact metric spaces. A map $f: X \ra Y$ is approximable by homeomorphisms (ABH) if, for every $\epsilon > 0$, $\exists$ a homeomorphism: $g: X \ra Y$ with $\operatorname{dist}_Y(f(x),g(x)) < \epsilon$ for all $x \in X$.

\medskip
\noindent
\textbf{Bing Shrinking Criterion (BSC)}: $f: X \ra Y$ is ABH iff for all $\epsilon > 0$ there is a homeomorphism $h: X \ra X$ so that
\begin{enumerate}
    \item $\operatorname{dist}_Y(f \circ h(x), f(x)) < \epsilon$ for all $x \in X$ and
    \item $\operatorname{diam}_X(h(f^{-1}(y))) < \epsilon$ for all $y \in Y$.
\end{enumerate}

Andrew and Rubin found self-homeomorphisms $h_\epsilon$ of $(D^2 \times S^1) \times \R$ so that $\pi \circ h_\epsilon$ is $\epsilon$-close to $\pi$, and $\operatorname{diam}(h(\pi^{-1}(r))) < \epsilon$, where $\pi: (D^2 \times S^1) \times \R \xrightarrow{\epsilon} (D^2 \times S^1 \slash \operatorname{Wh}) \times \R$ is the projection crushing each $\operatorname{Wh} \times r$ to a point.

Their idea is to lift, at any deep stage $i$, the self clasp of the Whitehead double $(D^2 \times S^1)_i \times r$ to an infinite chain of clasps and then apply a global twist making each chain of the link small diameter. Figure \ref{fig:lift} gives a hint of how this works.

\begin{figure}[ht]
    \centering
    \begin{tikzpicture}[scale=1.1]
        \draw[dashed] (-4,0) arc (0:180:1 and 0.25);
        \draw (-6,4) -- (-6,0) arc (180:360:1 and 0.25) -- (-4,4);
        \draw (-5,4) ellipse (1 and 0.25);
        \draw[dashed] (-4,1.5) arc(0:180:1 and 0.35);
        \draw[dashed] (-4.2,1.5) arc(0:180:0.8 and 0.15);
        \draw (-5.8,1.5) to[out=-80,in=180] (-5.3,1.3) to[out=0,in=160] (-5.1,1.3);
        \draw (-4.2,1.5) to[out=-100,in=0] (-4.95,1.33) to[out=180,in=155] (-5.05,1.15);
        \draw (-4.95,1.3) to[out=-45,in=0] (-5,1.1) to[out=180,in=-90] (-6,1.5);
        \draw(-4.9,1.1) to[out=0,in=-90] (-4,1.5);
    
        \draw[dashed] (-4,3) arc(0:180:1 and 0.35);
        \draw[dashed] (-4.2,3) arc(0:180:0.8 and 0.15);
        \draw (-5.8,3) to[out=-80,in=180] (-5.3,2.8) to[out=0,in=160] (-5.1,2.8);
        \draw (-4.2,3) to[out=-100,in=0] (-4.95,2.83) to[out=180,in=155] (-5.05,2.65);
        \draw (-4.95,2.8) to[out=-45,in=0] (-5,2.6) to[out=180,in=-90] (-6,3);
        \draw(-4.9,2.6) to[out=0,in=-90] (-4,3);
        \node at (-6.4,2) {$\R$};
        \draw[->] (-6.4,2.3) -- (-6.4,3.2);
        \draw[->] (-6.4,1.7) -- (-6.4,0.8);
        
        \draw[->] (-3.4,2) -- (-2.6,2);
        \node at (-3,2.2) {lift};
        \draw[->] (-3.7,4.3) -- (-4,4.1);
        \node at (-3,4.3) {$i$th stage};
        \node at (-3,3.3) {\scriptsize{$(i+1)$th stage,}};
        \node at (-3,2.8) {\scriptsize{$(D^2 \times S^1)_j \times h$}};
    
        \draw (-2,4) -- (-2,0) arc (180:360:1 and 0.25) -- (0,4);
        \draw (-1,4) ellipse (1 and 0.25);
        \draw[dashed] (0,0) arc (0:180:1 and 0.25);
        
        \foreach \y in {-1,-0.5,...,1} {
            \draw (-0.9,2.2+\y) to[out=195,in=170] (-1.05,2+\y) to[out=-10,in=190] (-0.95,2+\y);
            \draw (-0.9,2.2+\y) to[out=15,in=190] (0,2.35+\y);
            \draw (-1,2.075+\y) to[out=10,in=170] (-0.9,2.075+\y) to[out=0,in=15] (-1.05,1.9+\y) to[out=195,in=10] (-2,1.75+\y);
            \draw (-0.8,2.05+\y) to[out=10,in=190] (0,2.175+\y);
            \draw (-1.15,2.05+\y) to[out=190,in=10] (-2,1.925+\y);
        }
        \foreach \y in {-1.5,-1,...,0.5} {
            \draw[dashed] (0,2.35+\y) -- (-0.8,2.5+\y);
            \draw[dashed] (-1.2,2.62+\y) -- (-1.75,2.75+\y);
            \draw[dashed] (-0.5,2.32+\y) -- (-0.92,2.4+\y);
            \draw[dashed] (-1.37,2.55+\y) -- (-1.9,2.67+\y);
        }
        \draw[dashed] (-1.75,3.25) -- (-1.95,3.3);
        \draw[dashed] (-0.5,0.82) -- (-0.05,0.73);
        
        \draw[->] (0.4,2) -- (1.6,2);
        \node at (1,2.2) {twist};
    
        \draw (2,4) -- (2,0) arc (180:360:1 and 0.25) -- (4,4);
        \draw (3,4) ellipse (1 and 0.25);
        \draw[dashed] (4,0) arc (0:180:1 and 0.25);
        
        \foreach \x in {0,0.6,1.2} {
        \foreach \y in {-1.3,-0.8,...,1.2} {
            \draw (2.25+\x,2.4+\y) to[out=-95,in=180] (2.32+\x,1.97+\y) to[out=0,in=-90] (2.4+\x,2.1+\y);
            \draw (2.3+\x,2.1+\y) to[out=90,in=180] (2.38+\x,2.23+\y) to[out=0,in=85] (2.45+\x,1.8+\y);
        }
        \draw (2.465+\x,3.6) -- (2.45+\x,3.47);
        \draw (2.27+\x,0.65) -- (2.25+\x,0.52);
        }
        
        \node at (5.5,3) {$i$th stage solid};
        \node at (5.75,2.5) {tori appear small};
        \node at (5.55,2) {after lift/twist};
        \node at (5.95,1.5) {and contain $\operatorname{Wh} \times r$};
        \draw[->] (4.3,3) -- (3.75,2.8);
    \end{tikzpicture}
    \caption{}
    \label{fig:lift}
\end{figure}

The next step in the proof is what I called the design, but a better word would be a ``topographic map with some holes in it,'' perhaps singed in a campfire---but still useful. The map can be thought of as a certain closed subset $D$ of the standard (open) 2-hanlde $(D^2 \times \R^2,S^1 \times \R^2)$, but using the re-imbedding lemma, and its concomitant geometric control is also a subset of the (general Casson handle) $C$. Actually this topographic map varies in ways that are not terribly important with different choices of $C$, but in all cases the design $D$ includes into both $C$ and the standard (open) handle $H$. If this map was fully extensive and covered both $C$ and $H$ it would by itself describe the desired homeomorphism, but unfortunately it has a countable number of bits missing, which we call ``holes'' in $H$, and ``gaps'' in $C$. The gaps are bits of $C$ which remain terra incognita after the exploration that lead to the design. To sketch the rest of the proof, I need to tell you how the design is built, why it has these gaps, and finally what to do about them.

There will be a lot a branching going on in this paragraph, and I hope to keep two quite different types straight as they are described. Let me call them ``bulk'' and ``radial'' branching. In Figure \ref{fig:cassonhandle} we already met bulk branching, this comes from the fact that every Whitney disk that Casson installs to kill a double point has itself many double points and requires that many Whitney disks at the next level. Bulk branching, in the presence of geometric control, limits to a Cantor set (which plays the role of the singular points in the Andrews-Rubin argument). This branching is really only a distraction and does not materially affect the proof. The reader would be served well by the (false) fantasy that all Casson handles are unbranched (even when obtained via the re-imbedding theorem) and just forget about bulk branching. The more interesting branching is what I’m calling radial. Recall a 7-stage tower contains a 14-stage tower (Casson handles are alive and can replicate!). We can think of $14=7+7$ with a second 7-stage tower mounted on the top of the first. But the same combinatorial technology that allows raising height also allows us to build two alternative choices, $T_7$ and $T^\pr_7$ for the second tower with $T_7^\pr$ contained in $T_7$. This binary choice now occurs countably often as we may now identify 14 stage towers within both $T_7^\pr$ and $T_7$, and again create a binary choice (move in or stay out) in attaining the final 7 stages in these 14-ers. Continuing in this way we realize each of two choices ``in or out'' every 7th stage, always while maintaining geometric control so that the limits are as we expect. We call this branching ``radial'' because when $D$ is embedded in $H$, indeed at every such branch, the primed choice has a smaller $\rho$ coordinate in the polar coordinates of the transverse $\R^2$. The design carries a singular foliation. At radial coordinates which are not in the standard Cantor set the leaf is a compact 3-maniold with boundary, at a Cantor set radius the leaf is an Andrew-Rubin-like decomposition space. Corresponding to the ``middle thirds'' are the holes or gaps where this procedure has not succeeded in corresponding points of $H$ to points of $C$.

There is now a technical point, all but one of the holes is homeomorphic to $S^1 \times \R^3$ (the outlier is homeomorphic to $\R^4$). At this point we don’t know what the topology of the gaps is-–-but long after this proof is finished we do learn enough to know that they all along had the topology of their corresponding holes. Our plan for the $\{$holes$\}$ and the corresponding $\{$gaps$\}$ is eventually to crush them to points and analyze the results. It is much better in applying the Bing/Edwards mathematics that we crush only closed ``cell-like'' \cite{edwards80} sets. Cell-like means that when embedded in a Euclidean space (or Hilbert space) for any neighborhood $U$ there is a smaller neighborhood $V$ that contracts within $U$. Our holes/gaps are neither closed nor cell-like. The first deficiency is easily corrected by taking the closures, the second by locating certain 2-disks to add to the holes or gaps. We use the notations hole+ and gap+ for holes or gaps that have been closed and augmented with an additional 2-disk to become cell-like. In a sense adding the + is taking a step back because we are giving up some of the ``explored region'' of $C$ covered by the design, but sometimes to take two steps forward one must take a step backwards. This is an example.

Now comes the key diagram which will divide what remains of the proof into two separate theorems: 1) the Edwards shrink and 2) the ``Sphere to Sphere'' theorem.

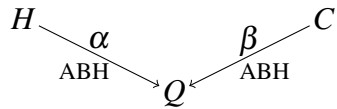
\begin{figure}[ht]
    \centering
    \begin{tikzpicture}
        \node at (0,0) {$Q$};
        \node at (2,1) {$C$};
        \node at (-2,1) {$H$};
        \draw[->] (-1.8,0.9) -- (-0.2,0.1);
        \node at (-1,0.7) {$\alpha$};
        \node at (-1.2,0.3) {\scriptsize{ABH}};
        \node at (1,0.7) {$\beta$};
        \node at (1.2,0.3) {\scriptsize{ABH}};
        \draw[->] (1.8,0.9) -- (0.2,0.1);
    \end{tikzpicture}
    \caption{$\alpha$ crushes $\{$holes+$\}$ and $\beta$ crushes $\{$gaps+$\}$.}
\end{figure}

Tautologically $H$ and $C$ have a common quotient $Q$, obtained by crushing the holes+, respectively gaps+, to points. This is the benefit of having found the design $D$ on both sides. It turns out both quotient maps can be shown to be approximable by homeomorphisms (ABH). This, of course, would finish the proof that Casson handles $C$ are homeomorphic to the standard open handle $H$, by composing one approximating homeomorphism by the inverse of the other. And with this result in hand the Poincar\'{e} Conjecture PC4, and much more within 4D topology follows.
\bigskip

Edwards analysis of the first map $\alpha$ is a tour de force of Bing topology. Many of the best ideas in that field enter: It encompasses and generalizes the Andrews-Rubin shrinking technology and also contains a compound use of the principle that ``Countable, null, star-like equivalent'' decompositions are shrinkable \cite{bean67}. It is a beautiful and pure application of the Bing shrinking criterion.

\begin{figure}[ht]
    \centering
    \begin{tikzpicture}[scale=1.4]
        \draw (-1,0) -- (6,0);
        \draw (0,-1) -- (0,3.5);
        
        \draw (1.9,-0.1) -- (4.2,-0.1) -- (3.7,-1.3) -- (1.4,-1.3) -- cycle;
        \draw (1.9,-0.7) circle (0.05);
        \draw[pattern=vertical lines] (2.05,-0.75) to[out=90,in=180] (2.5,-0.45) to[out=0,in=180] (2.9,-0.35) to[out=0,in=100] (3.3,-0.65) to[out=-80,in=0] (3.1,-0.85) to[out=180,in=0] (2.7,-0.95) to[out=180,in=-90] (2.05,-0.75);
        \draw (3.7,-0.6) circle  (0.05);
        \draw (3.45,-0.5) ellipse (0.1 and 0.07);
        \path[fill=white] (2.7,-0.65) circle (0.2);
        \node at (2.7,-0.65) {$X$};
        
        \draw (4.6,-0.3) ellipse (0.17 and 0.1);
        \draw (4.95,-0.2) circle (0.05);
        \draw (0.7,-0.3) ellipse (0.17 and 0.1);
        \draw (1.05,-0.2) circle (0.05);
        
        \draw (0,0) -- (0.53,0.53) -- (0.87,0.53) -- (1,0.67) -- (1.1,0.67) -- (1.85,1.4) -- (1.95,1.4) -- (2.05,1.5) -- (3.3,1.5) -- (3.35,1.55) -- (3.55,1.55) -- (3.65,1.65) -- (3.75,1.65) -- (4.43,2.35) -- (4.77,2.35) -- (4.9,2.5) -- (5,2.5) -- (5.7,3.2);
        
        \draw[dashed] (0.53,0.53) -- (0.53,-0.3);
        \draw[dashed] (0.87,0.53) -- (0.85,-0.3);
        \draw[dashed] (1,0.67) -- (1,-0.2);
        \draw[dashed] (1.1,0.67) -- (1.1,-0.2);
        \draw[dashed] (1.85,1.4) -- (1.85,-0.7);
        \draw[dashed] (1.95,1.4) -- (1.95,-0.7);
        \draw[dashed] (2.05,1.5) -- (2.05,-0.75);
        \draw[dashed] (3.3,1.5) -- (3.3,-0.65);
        \draw[dashed] (3.35,1.55) -- (3.35,-0.5);
        \draw[dashed] (3.55,1.55) -- (3.55,-0.5);
        \draw[dashed] (3.65,1.65) -- (3.65,-0.5);
        \draw[dashed] (3.75,1.65) -- (3.75,-0.5);
        \draw[dashed] (4.43,2.35) -- (4.43,-0.3);
        \draw[dashed] (4.77,2.35) -- (4.77,-0.3);
        \draw[dashed] (4.9,2.5) -- (4.9,-0.2);
        \draw[dashed] (5,2.5) -- (5,-0.2);
        
        \draw[rotate=-90,xshift=-4.5cm, pattern=horizontal lines] (2.05,-0.75) to[out=90,in=180] (2.5,-0.45) to[out=0,in=180] (2.9,-0.35) to[out=0,in=100] (3.3,-0.65) to[out=-80,in=0] (3.1,-0.85) to[out=180,in=0] (2.7,-0.95) to[out=180,in=-90] (2.05,-0.75);
        \path[fill=white] (-0.67,1.8) circle (0.2);
        \node at (-0.67,1.8) {$X$};
        \draw[rotate=-90,xshift=-4.5cm] (1.9,-0.1) -- (4.2,-0.1) -- (3.7,-1.3) -- (1.4,-1.3) -- cycle;
        \draw[rotate=-90,xshift=-4.5cm]  (1.9,-0.7) circle (0.05);
        \draw[rotate=-90,xshift=-4.5cm] (3.7,-0.6) circle  (0.05);
        \draw[rotate=-90,xshift=-4.5cm] (3.45,-0.5) ellipse (0.1 and 0.07);
        
        \node at (-0.7,3.8) {$R^4 \supseteq Q \cong H$};
        \node at (3,2.5) {flat spot$(s)$};
        \draw[->] (3,2.3) -- (2.7,1.6);
        \node at (2.75,-1.7) {graph of $\beta$};
        \draw[->] (1.25,-1.1) to[out=180,in=235] (-1.5,1);
        \node at (6,-0.8) {$C =$ Casson handle};
    \end{tikzpicture}
    \caption{The neighborhood of the largest gap+, $X$, is drawn on the $y$-axis near $\beta(X)$ to facilitate a resolution of $\beta$ to a homeomorphism on $X$. This is good, but small vertical spots inevitably result from pulling points where $\beta$ is already a homeomorphism over other smaller elements of $\{$gaps+$\}$. This deficit is corrected later.}
    \label{fig:procedure}
\end{figure}
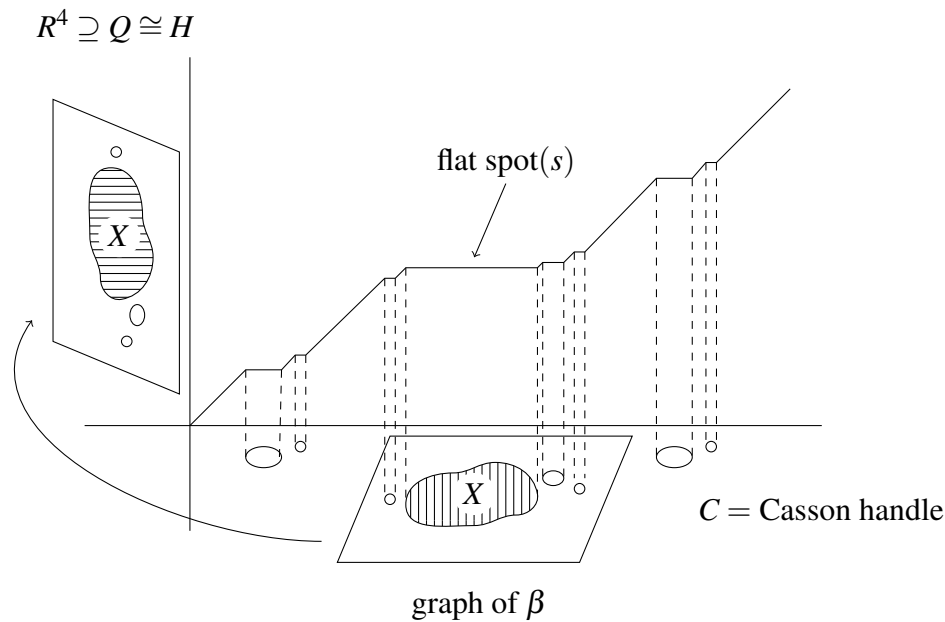

The arrow $\beta$, might be called the ``blindfold shrink.''  Ric Ancel \cite{ancel84}, Jim Cannon, Frank Quinn, and Mike Starbird helped me revise this argument and express it as an extension of Brown’s proof of the Schoenflies theorem, a proof that I learned from them. The notation is a bit daunting but the idea is simple. After the Edwards shrink we know that both domain and range of $\beta$ are manifolds and ones that embed in $\R^4$, in fact one readily reduces to the case that both domain and range are standard 4-spheres. Then the idea is to look at the 8D graph of the function $\beta$ and try to improve it to a homeomorphism by systematically modifying it to remove its largest ``flat spots.'' The data (from the Edwards shrink) that the quotient space has nice local and global topology and allows you to insert a local ``drawing'' of the decomposition space structure into the target within a small neighborhood of the image of these largest decomposition elements. A natural procedure, see Figure \ref{fig:procedure}, uses this inserted drawing to resolve these most singular points and make the function one-to-one over them. There is a blending problem: the resolution creates countably many small vertical spots, that is the resolution is only a relation not a function. But one should not worry, we are on the right track. One goes back and forth: sanding down first the largest horizontal spots, then flipping the workpiece over and sanding down the largest of these newly created vertical spots, flipping it back to sand down somewhat smaller horizontal spots, etc. One moves back and forth sanding and polishing until the limiting relation has neither horizontal nor vertical thickness, i.e.\ becomes a homeomorphism.

I thank Arumina Ray and Mark Powell for including this lecture in their book; it is hoped that these brief recollections may add context or at least amusement. There is much left to be understood in four dimensions---the best of luck!

\newpage
\bibliography{references}

\begin{bibdiv}
\begin{biblist}

\bib{ancel84}{article}{
      author={Ancel, Fredric},
       title={Approximating cell-like maps of ${S}^4$ by homeomorphisms},
        date={1984},
     journal={Four-manifold Theory (Durham, NH, 1982)},
      volume={35},
       pages={143\ndash 164},
}

\bib{andrews65}{article}{
      author={Andrews, J.J.},
      author={Rubin, Leonard},
       title={Some spaces whose product with ${E}^1$ is ${E}^4$},
        date={1965},
     journal={Bull. Amer. Math. Soc.},
      volume={71},
      number={4},
       pages={675\ndash 677},
}

\bib{bean67}{article}{
      author={Bean, Ralph},
       title={Decompositions of ${E}^3$ with a null sequence of starlike
  equivalent nondegenerate elements are ${E}^3$},
        date={1967},
     journal={Illinois J. Math.},
      volume={11},
       pages={21\ndash 23},
}

\bib{casson86}{article}{
      author={Casson, Andrew},
       title={Three lectures on new infinite constructions in 4-dimensional
  manifolds},
        date={1986},
     journal={A la Recherche de la Topologie Perdue},
       pages={201\ndash 244},
}

\bib{cha16}{article}{
      author={Cha, Jae~Choon},
      author={Powell, Mark},
       title={Casson towers and slice links},
        date={2016},
     journal={Invent. Math.},
      volume={205},
      number={2},
       pages={413\ndash 457},
}

\bib{donaldson86}{article}{
      author={Donaldson, Simon},
       title={Connections, cohomology and the intersection forms of
  4-manifolds},
        date={1986},
     journal={J. Differential Geom.},
      volume={24},
      number={3},
       pages={275\ndash 341},
}

\bib{edwards80}{inproceedings}{
      author={Edwards, Robert},
       title={The topology of manifolds and cell-like maps},
        date={1980},
   booktitle={Proceedings of the {I}nternational {C}ongress of {M}athematicians
  ({H}elsinki, 1978)},
       pages={111\ndash 127},
}

\bib{fox66}{article}{
      author={Fox, Ralph},
      author={Milnor, John},
       title={Singularities of 2-spheres in 4-space and cobordism of knots},
        date={1966},
     journal={Osaka J. Math.},
      volume={3},
      number={2},
       pages={257\ndash 267},
}

\bib{gs84}{incollection}{
      author={Gompf, Robert},
      author={Singh, Sukhjit},
       title={On {F}reedman's reimbedding theorems},
        date={1984},
   booktitle={Four-manifold theory},
      editor={Gordon, Cameron},
      editor={Kirby, Robion},
      series={Contemporary Mathematics},
      volume={35},
   publisher={American Mathematical Soc.},
}

\bib{lashof70}{inproceedings}{
      author={Lashof, Richard},
       title={The immersion approach to triangulation and smoothing},
        date={1971},
   booktitle={Algebraic topology ({P}roc. {S}ympos. {P}ure {M}ath., {V}ol.
  {XXII}, {U}niv. {W}isconsin, {M}adison, {W}is., 1970)},
       pages={131\ndash 164},
}

\bib{moise52a}{article}{
      author={Moise, Edwin~E.},
       title={Affine structures in 3-manifolds: {IV}. {P}iecewise linear
  approximations of homeomorphisms},
        date={1952},
     journal={Ann. of Math.},
      volume={55},
      number={2},
       pages={215\ndash 222},
}

\bib{moise52b}{article}{
      author={Moise, Edwin~E.},
       title={Affine structures in 3-manifolds: {V}. {T}he triangulation
  theorem and {H}auptvermutung},
        date={1952},
     journal={Ann. of Math.},
      volume={56},
      number={1},
       pages={96\ndash 114},
}

\bib{perelman03a}{article}{
      author={Perelman, Grisha},
       title={Finite extinction time for the solutions to the ricci flow on
  certain three-manifolds},
        date={2003},
      eprint={arXiv:math/0307245},
}

\bib{perelman03b}{article}{
      author={Perelman, Grisha},
       title={Ricci flow with surgery on three-manifolds},
        date={2003},
      eprint={arXiv:math/0303109},
}

\bib{smale56}{article}{
      author={Smale, Stephen},
       title={Generalized {P}oincar{\'e} conjecture in dimensions greater than
  4},
        date={1956},
     journal={Ann. of Math.},
      volume={64},
       pages={399\ndash 405},
}

\bib{whitney44}{article}{
      author={Whitney, Hassler},
       title={The self-intersections of a smooth n-manifold in 2n-space},
        date={1944},
     journal={Ann. of Math.},
      volume={45},
      number={2},
       pages={220\ndash 246},
}

\end{biblist}
\end{bibdiv}

\end{document}